\documentclass[reqno,12pt]{amsart}

\usepackage{amsmath, amssymb, amsthm, amsfonts} 
\usepackage[english]{babel}
\usepackage{bbm}
\usepackage{graphicx}
\usepackage{url}
\usepackage{epstopdf}
\usepackage[a4paper,bindingoffset=0.5cm,left=2cm,right=2cm,top=2.5cm,bottom=2cm,footskip=.8cm]{geometry}

\usepackage{tikz}
\usetikzlibrary{arrows, automata,positioning,calc,shapes,decorations.pathreplacing,decorations.markings,shapes.misc,petri,topaths}
  \usepackage{pgfplots}
  \pgfplotsset{compat=newest}
  \usetikzlibrary{plotmarks}
  \usepackage{grffile}
\newlength\figureheight
  \newlength\figurewidth
\pgfplotsset{%
    tick label style={font=\scriptsize},
    label style={font=\footnotesize},
    legend style={font=\footnotesize},
         every axis plot/.append style={very thick}
}

\usepackage{rotating}
\usepackage{amsbsy,enumerate}
\usepackage{graphicx}
\usepackage{ccaption}
\usepackage{comment}
\usepackage{mathrsfs} 
\usepackage{diagbox}
\usepackage{caption}

\newcommand{\vb}{\vspace{3.2mm}}
\renewcommand{\hat}{\widehat}

\allowdisplaybreaks

\newtheorem{theorem}{Theorem}
\newtheorem{remark}{Remark}

\newtheorem{proposition}{Proposition}

\makeatletter
\renewcommand{\fnum@figure}[1]{\textbf{\figurename~\thefigure}. }
\renewcommand{\fnum@table}[1]{\textbf{\tablename~\thetable}. }
\makeatother

\begin{document}

\title[L\'evy processes inspected at Poisson moments]{A decomposition for L\'evy processes\\ inspected at Poisson moments}

\author{Onno Boxma {and} Michel Mandjes}

\begin{abstract}
We consider a L\'evy process $Y(t)$ that is not permanently observed, but rather inspected at Poisson($\omega$) moments only, over an exponentially distributed time $T_\beta$ with parameter $\beta$. The focus lies on the analysis of the distribution of the running maximum at such inspection moments up to $T_\beta$, denoted by $Y_{\beta,\omega}$.
Our main result is a decomposition:  we derive a remarkable distributional equality that contains 
$Y_{\beta,\omega}$ as well as the running maximum process $\bar Y(t)$ at 
the exponentially distributed times $T_\beta$ and $T_{\beta+\omega}$.
Concretely, 
$\overline{Y}(T_\beta)$ can be written the sum of the two independent random variables that are distributed as $Y_{\beta,\omega}$ and $\overline{Y}(T_{\beta+\omega})$.
The distribution of $Y_{\beta,\omega}$ can be identified more explicitly in the two special cases of a spectrally positive and a spectrally negative L\'evy process. As an illustrative example of the potential of our results, we show how to determine the asymptotic behavior of the bankruptcy probability in the Cram\'er-Lundberg insurance risk model.

\vb


\vb

\noindent
{\sc Affiliations.} 
O.\ Boxma is with E{\sc urandom} and the Department of Mathematics and Computer Science, Eindhoven University of Technology, Eindhoven, the Netherlands.  
M.\ Mandjes is with the Korteweg-de Vries Institute for Mathematics, University of Amsterdam, Science Park 904, 1098 XH Amsterdam, the Netherlands; he is also with E{\sc urandom}, Eindhoven University of Technology, Eindhoven, the Netherlands, and Amsterdam Business School, Faculty of Economics and Business, University of Amsterdam, Amsterdam, the Netherlands.

\noindent
This research is partly funded by the NWO Gravitation project N{\sc etworks}, grant number 024.002.003. Email addresses: \url{o.j.boxma@tue.nl} and \url{m.r.h.mandjes@uva.nl}.

\vb

\noindent
{\sc Acknowledgments.} The authors thank Hansjoerg Albrecher and Jevgenijs Ivanovs for various useful comments.

\vb

\noindent
{\it Date}: {\today}.

\end{abstract}

\maketitle 
\newpage
\section{Introduction}
\noindent
We consider a general real-valued L\'evy process $Y \equiv \{Y(t), ~t \geqslant 0\}$ that is inspected at the epochs $I_1,I_2,\dots$ of an independently evolving Poisson process with intensity $\omega>0$.
Our aim is to compare the distributions of the running maximum of the L\'evy process under continuous observation and its counterpart at the inspection epochs. We do this until an, independently sampled, exp($\beta$) distributed `killing time' $T_\beta$ (for $\beta\geqslant 0$). This also covers the case of an infinite time horizon: let $\beta$ go to $0$, imposing the additional assumption of the L\'evy process' mean to be negative to avoid the running maximum drifting to $\infty$.

The motivation behind our study lies in the fact that in many real-life situations the stochastic process under study is not, or {\it can} not, be continuously observed,
but is rather inspected at discrete times. Examples abound in 
reliability and healthcare, where an object or person is checked regularly during the lifetime. Our interest in this topic mainly stems from {an application in} insurance risk, and in particular from the classical Cram\'er-Lundberg model. That model represents the surplus process of an insurance company that earns money (at a constant rate) via premiums, and that pays claims to its clients according to a compound Poisson process.
The main quantity of interest is the ruin probability $p(u)$, viz., the probability that the surplus level becomes negative when starting with initial capital $u>0$.
Albrecher, Gerber and Shiu \cite{AGS} relaxed this model, observing that companies sometimes can continue doing business even when they are technically ruined.
More specifically, they distinguished ruin and bankruptcy, the latter occurring when the surplus is negative {\em at an inspection epoch}. It is obvious that the probability of bankruptcy $\tilde{p}(u)$ is smaller than $p(u)$, but one would like to quantify the difference.
Albrecher and Lautscham \cite{AL} have determined the bankruptcy probability in the Cram\'er-Lundberg setting with exponentially distributed claim sizes;
in~\cite{BEJ} this was generalized to the case of generally distributed claim sizes. In both papers, the inspection rate was allowed to depend on the current surplus level.
In~\cite{BEJ}, and also in \cite{ABEK}, a related queueing (or inventory) model was also studied, where a server works even when there are no customers (or orders), building up storage that is removed at the Poisson inspection epochs.

While the results in the present paper allow one to get a better insight into the relation between the ruin and bankruptcy probability in the Cram\'er-Lundberg model
and some of its (L\'evy) generalizations, our motivation is also to a large extent theoretical.
We prove that there is a remarkably simple relation between the running maximum $\overline{Y}(T_{\beta})$ of the supremum of the L\'evy process $Y$  until the exp($\beta$) killing epoch
and the running maximum $Y_{\beta,\omega}$ of $Y$ at inspection epochs until that killing: our main result (Theorem~\ref{the2}) is that
\begin{equation}\label{doco}\overline Y({T_\beta}) \stackrel{\rm d}{=} Y_{\beta,\omega}+ \overline Y({T_{\beta+\omega}}),\end{equation} the two quantities in the righthand side being independent. The proof of this result relies on known results from Wiener-Hopf theory for L\'evy processes, and is surprisingly straightforward.

Our main result, i.e., the decomposition \eqref{doco}, is proven in Section~\ref{sec2}. In that section we also show how it can, alternatively, be obtained by applying
results from \cite{AI}. The latter paper also served as an important source of inspiration for us, as it presents several beautiful identities relating exit problems for L\'evy processes
under permanent observation and their counterparts under Poisson inspections.
While Section~\ref{sec2} focuses on a general L\'evy process, we restrict ourselves in Section~\ref{sec3} to spectrally one-sided L\'evy processes. In the spectrally positive case, we succeed in expressing the transform of $Y_{\beta,\omega}$ explicitly in terms of the Laplace exponent of the driving L\'evy process $Y$. In the spectrally negative case we prove that $Y_{\beta,\omega}$ is exponentially distributed, with an atom at zero.
As an illustration of the application potential of our results, we study in Section~\ref{sec4} the asymptotics (for large initial capital $u$, that is) of the bankruptcy probability in the setting of the Cram\'er-Lundberg model, distinguishing between the cases of light-tailed and heavy-tailed claim sizes. At various places we explicitly exploit the relation between fluctuation-theoretic concepts and their queueing counterparts, effectively resting on a duality between L\'evy-type insurance risk models and corresponding queueing models; regarding this duality, see e.g.\ the account in \cite[Section III.2]{AA}.
In addition, we frequently use the fact that, due to the Wiener-Hopf decomposition, the increment of $Y$ between two inspection epochs can be written as the difference between two independent positive random variables, which enables us to write the quantities under study in terms of waiting times in associated queueing models.

\section{Decomposition}
\label{sec2}
This section establishes the decomposition \eqref{doco}, and relates it to results that have recently appeared in \cite{AI}.

\subsection{The main decomposition result}
\label{sec2.1}

\vb 

\noindent
Let $Y\equiv \{Y(t), ~t \geqslant 0\} $ be a general real-valued L\'evy process. In addition, let $\bar Y$ be the associated running maximum process:
\begin{equation}\overline Y(t) := \sup_{s \in [0,t]} Y(s).\end{equation} 
We consider the setting in which, for some intensity $\omega>0$, the L\'evy process is inspected at Poisson($\omega$) moments $I_1,I_2,\dots$, so that the number of inspections up to time $t$ has a Poisson distribution with mean $\omega t$. We consider the resulting inspected process until `killing', which happens at an exponentially distributed time $T_\beta$ with parameter $\beta\geqslant 0$, sampled independently from $Y$.

To analyze the process $Y$ at inspection moments, we denote by $Z_m$ the increment of the L\'evy process $Y(t)$ between the two consecutive inspection times $I_{m-1}$ and $I_m$ (with $I_0 :=0$), conditioned on the process not having been killed. 
We also define $S_0:=0$ and,
for $n\in{\mathbb N}$,
\begin{equation}S_n:=\sum_{m=1}^n Z_m,\end{equation}
and the corresponding running maximum process $\overline S_n:=\max\{S_0,S_1,\ldots,S_n\}.$
We wish to analyze the running maximum of the inspected process until killing. The number of inspections $N_{\beta,\omega}$ before killing is shifted-geometric, with the `killing probability' given by $\beta/(\beta+\omega)$:
\begin{equation}
{\mathbb P}(N_{\beta,\omega} = n) = \left(\frac{\omega}{\beta+\omega}\right)^n \frac{\beta}{\beta+\omega}.
\end{equation}
The random variable of our interest is
\begin{equation}
Y_{\beta,\omega} := \overline S_{N_{\beta,\omega}} = \sup_{n=0,1,\ldots,N_{\beta,\omega}} S_n.
\label{2.1}
\end{equation}
Observe that we have the identity
\begin{equation}{\mathbb E}\,{\rm e}^{-\alpha Y_{\beta,\omega}} = \sum_{n=0}^\infty 
\left(\frac{\omega}{\beta+\omega}\right)^n \frac{\beta}{\beta+\omega} \,{\mathbb E}\,{\rm e}^{-\alpha \overline S_n}.\end{equation}

\begin{theorem}
\label{the1}
For any $\alpha>0$,
\begin{equation}{\mathbb E}\,{\rm e}^{-\alpha Y_{\beta,\omega}} = \exp\left(-\int_0^\infty\int_{(0,\infty)} \frac{1}{t} {\rm e}^{-\beta t}(1-{\rm e}^{-\omega t}) (1-{\rm e}^{-\alpha x})\, \,{\mathbb P}(Y(t)\in {\rm d}x)\,{\rm d}t
\right).
\label{eq7}
\end{equation}
\end{theorem}

\noindent {\it Proof.}
As pointed out in \cite{Kyprianou2}, applying Wiener-Hopf theory for random walks,
\begin{equation}\sum_{n=0}^\infty  
(1-p)^n p \,{\mathbb E}\,{\rm e}^{-\alpha \overline S_n} = 
\exp\left(-\int_{(0,\infty)} \sum_{n=1}^\infty \frac{1}{n} (1-{\rm e}^{-\alpha x})(1-p)^n\,{\mathbb P}(S_n\in {\rm d}x)\right);\end{equation}
see also e.g.\ \cite[Section 3.3]{DM}. As a consequence, 
\begin{equation}\label{trafmax}{\mathbb E}\,{\rm e}^{-\alpha Y_{\beta,\omega}} = \exp\left(-\int_{(0,\infty)} \sum_{n=1}^\infty \frac{1}{n} (1-{\rm e}^{-\alpha x})\left(\frac{\omega}{\beta+\omega}\right)^n\,{\mathbb P}(S_n\in {\rm d}x)\right).\end{equation}
Now realize that, with ${\rm E}(n,a)$ an Erlang random variable with shape parameter $n\in{\mathbb N}$ and scale parameter $a>0$, conditional on the process not having been killed,
\begin{equation}S_n\stackrel{\rm d}{=} Y({\rm E}(n, \beta+\omega)),\end{equation}
so that
\begin{equation} {\mathbb P}(S_n\in {\rm d}x) = \int_0^\infty (\beta+\omega)^n t^{n-1} \frac{e^{-(\beta+\omega)t}}{(n-1)!}\,{\mathbb P}(Y(t)\in {\rm d}x)\,{\rm d}t.\end{equation}
We thus obtain
\begin{align}{\mathbb E}\,{\rm e}^{-\alpha Y_{\beta,\omega}} &= \exp\left(-\int_0^\infty\int_{(0,\infty)} \frac{1}{t}\sum_{n=1}^\infty \frac{1}{n!} (1-{\rm e}^{-\alpha x})({\omega}t)^n\,{\rm e}^{-(\beta+\omega)t}\,{\mathbb P}(Y(t)\in {\rm d}x)\, {\rm d}t
\right)
\nonumber
\\
&= \exp\left(-\int_0^\infty\int_{(0,\infty)} \frac{1}{t}({\rm e}^{\omega t}-1) (1-{\rm e}^{-\alpha x})\,{\rm e}^{-(\beta+\omega)t} \,{\mathbb P}(Y(t)\in {\rm d}x)\,{\rm d}t
\right)
\nonumber
\\
&= \exp\left(-\int_0^\infty\int_{(0,\infty)} \frac{1}{t} {\rm e}^{-\beta t}(1-{\rm e}^{-\omega t}) (1-{\rm e}^{-\alpha x})\, \,{\mathbb P}(Y(t)\in {\rm d}x)\,{\rm d}t
\right).
\end{align}
This proves the claim. \hfill$\Box$

\vb

We now state and prove our main result, a decomposition theorem for a L\'evy process with Poisson inspection epochs.

\begin{theorem}
\label{the2}
The following distributional equality applies:
\begin{equation}
    \label{doco2}
\overline Y({T_\beta}) \stackrel{\rm d}{=} Y_{\beta,\omega}+ \overline Y({T_{\beta+\omega}}),\end{equation}
with the two terms in the righthand side being independent. 
\end{theorem}

\noindent {\it Proof.} By applying the Wiener-Hopf theory for L\'evy processes, as presented in e.g.\ \cite[Theorem 6.15]{Kyprianou} or \cite[Section 3.3]{DM}, for any $\zeta>0$,
\begin{equation}{\mathbb E} \,{\rm e}^{-\alpha \overline Y({T_\zeta})}= \exp\left(-\int_0^\infty\int_0^\infty \frac{1}{t} {\rm e}^{-\zeta t} (1-{\rm e}^{-\alpha x})\, \,{\mathbb P}(Y(t)\in {\rm d}x)\,{\rm d}t
\right).\end{equation}
Taking $\zeta=\beta$ and $\zeta = \beta + \omega$, and using Theorem~\ref{the1}, we obtain
\begin{equation}{\mathbb E} \,{\rm e}^{-\alpha \overline Y({T_\beta})} = {\mathbb E}\,{\rm e}^{-\alpha Y_{\beta,\omega}}\cdot {\mathbb E} \,{\rm e}^{-\alpha \overline Y({T_{\beta+\omega}})},\end{equation}
which implies the stated.\hfill$\Box$

\begin{remark}{\em 
 A striking aspect of the decomposition, besides its remarkably straightforward proof, is that the impact of $\omega$ in the first term in the righthand side of \eqref{doco2} apparently equals the impact of $\omega$ in the second term, but `with opposite sign' -- to this end, observe that the lefthand side of \eqref{doco2} does not
involve $\omega$ at all.
Observe that the first term in the righthand side is increasing in $\omega$ (as the inspection process takes place at an increasingly high frequency, with the length of the interval held fixed),
whereas the second term is decreasing in $\omega$ (as a supremum over an increasingly small interval is taken).} \hfill$\Diamond$
\end{remark}

\begin{remark}{\em 
As is to be expected,
$S_{N_{\beta,\omega}}$ can also be decomposed as the sum of $\overline S_{N_{\beta,\omega}}$ and $\underline S_{N_{\beta,\omega}}
:= \inf_{n=0,1,\ldots,N_{\beta,\omega}} S_n$,
with the latter two quantities being independent.
To verify this, first observe that
${\mathbb E} \, {\rm e}^{\alpha {\rm i} Z_m} = (\beta + \omega)/(\beta + \omega - {\rm log} \, {\mathbb E} \, {\rm e}^{\alpha{\rm i} Y(1)})$. 
Hence,
\begin{align}
   {\mathbb E}\,e^{\alpha{\rm i} S_{N_{\beta,\omega}}} &= \frac{\beta}{\beta+\omega -\omega\, {\mathbb E}\,e^{\alpha{\rm i}Z}} \, = \,
\frac{1 + \omega/(\beta - {\rm log} \, {\mathbb E} \, {\rm e}^{\alpha{\rm i} Y(1)})}{1 + \omega/\beta}
\nonumber
\\
    &= \exp\left(-\int_0^\infty\int_{(-\infty,\infty)} \frac{1}{t} {\rm e}^{-\beta t}(1-{\rm e}^{-\omega t}) (1-{\rm e}^{\alpha{\rm i} x})\, \,{\mathbb P}(Y(t)\in {\rm d}x)\,{\rm d}t
\right).
\end{align}
In the last step we have used the Frullani integral \cite[Lemma 1.7]{Kyprianou}
\begin{equation}
1 + \frac{\omega}{a} \, = \, {\rm exp} \left(\int_0^{\infty} \frac{1}{t} {\rm e}^{-at} (1 - {\rm e}^{-\omega t})  \, {\rm d}t \right),
\end{equation}
with both $a = \beta$ and $a = \beta - \log \, {\mathbb E}\,e^{\alpha{\rm i} Y(1)}$.
Finally observe that, by symmetry, the running minimum 
can be dealt with in the precise same manner as the running maximum, in that the transform of $\underline S_{N_{\beta,\omega}}$ is as given in \eqref{eq7}, but with the integration interval $(0,\infty)$ replaced by $(-\infty,0)$. 
}
\hfill$\Diamond$
\end{remark}

\subsection{Relation to a result of Albrecher and Ivanovs}
\label{sec2.2}

\noindent
In this subsection we outline the relation between Theorem 2 and some results of Albrecher and Ivanovs \cite{AI}.
In preparation of this, we take a closer look at the increments $Z_m$ between inspection epochs, and we mention the powerful concept of Wiener-Hopf factorization for L\'evy processes; see, e.g.,   \cite[Chapter 6]{Kyprianou}
or   \cite[Section 3.3]{DM}. 

The Wiener-Hopf decomposition  entails that each $Z_m$ can be written as $Z_m = Z_m^+ - Z_m^-$ with $Z_m^+$ and $Z_m^-$ independent and both non-negative.
Here $Z_m^+$ (resp.\ $Z_m^-$) is distributed as the supremum (resp.\ minus the infimum) of the L\'evy process $Y$, when started anew at zero at inspection epoch $I_{m-1}$, over the interval between $I_{m-1}$ and $I_m$ (whose length is exp(${\beta+\omega}$)).
To get some feeling for this property, it is useful to observe that a time-reversibility argument for L\'evy processes implies that, with $\underline Y(t)$ denoting the running minimum process, we have
\begin{equation}
Y(t) - \underline Y(t) = Y(t) - \inf_{s \in [0,t]}Y(s) = \sup_{s \in [0,t]} (Y(t) - Y(s)) \stackrel{{\rm d}}{=} \sup_{s \in [0,t]} Y(s) = \overline{Y}(t) 
\label{WH}
\end{equation}
(but, evidently, for a {\it given} $t$, $\underline Y(t)$ and $\overline{Y}(t)$ are {\it not} independent). In the sequel $Z^+$ (resp.\ $Z^-$) is a generic random variable distributed as $Z_m^+$ (resp.\ $Z_m^-$).

Albrecher and Ivanovs \cite{AI} consider a L\'evy process $X\equiv \{X(t), t\geqslant 0\}$, starting at $u$, which is also being inspected at Poisson($\omega$) epochs $I_0=0,I_1,\dots$.
If $X$ attains a negative value, then ruin is said to occur, whereas if it is negative at an inspection epoch, then bankruptcy is said to occur.
Recall that the all-time ruin probability starting at surplus level $u$ is denoted by $p(u)$, and the (obviously smaller) corresponding bankruptcy probability by $\tilde{p}(u)$.
The starting-point in \cite{AI} is their elegant Proposition 1, which (in our notation) states that
\begin{equation}
p(u) = {\mathbb E} \,\tilde{p}(u-Z^+) ,
\label{d1}
\end{equation}
where the process $X$ relates to our $Y$ through the relation $X(t)=u-Y(t)$ for $t\geqslant 0$.
The focus in \cite{AI} lies not on deriving decompositions, but the proof of (\ref{d1}) in fact implicitly reveals such a decomposition, which can be used so as to rederive our decomposition of Theorem \ref{the2}.
They introduce (again adapted to our notation) the partial sums, for $i=1,2,\dots$,
\begin{align}
& \hat{\sigma}_0 := 0, \,\,\,\,\:\: \hat{\sigma}_i := - \sum_{j=0}^i Z_j,
\label{d2}
\\
& \sigma_0 := -Z_1^+, \,\,\,\,\:\: \sigma_i := - \sum_{j=1}^i (Z_{j+1}^+ -Z_j^-) .
\label{d3}
\end{align}
Then it is concluded that
\begin{equation}
\{\hat{\sigma}_i -Z^+ \}_{i=0,1,\ldots} \stackrel{{\rm d}}{=} \{\sigma_i\}_{i=0,1,\ldots} ,
\label{d4}
\end{equation}
and hence 
\begin{equation}
p(u) = {\mathbb P}\big(-  \min_{i \geqslant 0} \, \sigma_i \geqslant u\big) = {\mathbb P}\big(- \min_{i \geqslant 0} \, \hat{\sigma}_i + Z_1^+ \geqslant u\big) = {\mathbb E} \,\tilde{p}(u-Z_1^+) .
\label{d5}
\end{equation}
Observe that this identity implicitly entails the decomposition
\begin{equation}- \min_{i=0,1,\ldots} \, \sigma_i \,  \stackrel{{\rm d}}{=} \, - \min_{i =0,1,\ldots} \,  \hat{\sigma}_i + Z_1^+ .\end{equation}

Let us now turn to the  case with `killing', as has been considered in Theorem \ref{the2}, i.e., the process ends at $T_\beta \sim {\rm exp}(\beta)$.
We note that \cite[Remark 3]{AI} briefly mentions the option of killing, at an inspection epoch $I_i$. It is stated there, without proof, that the finite-time
ruin and bankruptcy probabilities before inspection time $I_i$ are related via
\begin{equation}
p(u,I_i) = {\mathbb P}\big(- \min_{j=0,1,\ldots,  i-1}  \sigma_j \geqslant u\big) = {\mathbb P}\big(- \min_{j=0,1,\ldots,  i-1} \hat{\sigma}_j + Z_1^+ \geqslant u\big) = {\mathbb E} \,\tilde{p}(u-Z_1^+,I_{i-1}) ;
\label{d6}
\end{equation}
here $p(u,t)$ (resp.\ $\tilde p(u,t)$) is the probability of ruin (resp.\ bankruptcy) before time $t$, given an initial surplus $u$.
To translate this observation to the setting of Theorem \ref{the2}, let us assume that the system is inspected at the Poisson($\beta+\omega$) epochs $I_0=0,I_1,\dots$.
The inspection intervals now are exp($\beta+\omega$) distributed, and accordingly, in the distributions of the $Z^+$ and $Z^-$  the parameter $\omega$ should be replaced by $\beta + \omega$.
The `$\beta$-inspection' is preceded by $N_{\beta,\omega}$ `$\omega$-inspections' at epochs $I_1,\dots,I_{N_{\beta,\omega}}$.
Now consider the three terms in Theorem~\ref{the2}, and compare them with the three main random elements featuring in (\ref{d6}).
Firstly, observe that
\begin{equation}
Z_1^+ \stackrel{{\rm d}}{=} \overline{Y}(T_{\beta+\omega}),
\label{d6a}
\end{equation}
as both are distributed as the supremum of the L\'evy process $Y$ over an exp($\beta + \omega$) interval.
Secondly, noticing that if $I_i = T_\beta$ then $N_{\beta,\omega} = i-1$, we have
\begin{equation}
- \min_{j=0,1, \ldots, N_{\beta,\omega}} \, \hat{\sigma}_j = \sup_{j=0,1, \ldots, N_{\beta,\omega}} \sum_{k=0}^j Z_k \stackrel{{\rm d}}{=} \overline{Y}_{\beta,\omega} .
\label{d7}
\end{equation}
Thirdly, 
\begin{equation}
- \min_{j=0,1, \ldots, N_{\beta,\omega}} \, \sigma_j = \sup_{j=0,1, \ldots, N_{\beta,\omega}} \sum_{k=0}^j (Z_{k+1}^+ - Z_k^-) \stackrel{{\rm d}}{=} \overline{Y}(T_\beta) ,
\label{d8}
\end{equation}
as the latter supremum is the supremum of the $Y$ process until $T_\beta$. We thus conclude that Theorem \ref{the2} can in this way be recovered from the middle equality in (\ref{d6}).

\section{The two spectrally one-sided cases}
\label{sec3}
In this section we consider two special cases for which  the various components of Theorem~\ref{the2} can be obtained through an explicit characterization.
Section~\ref{sec3.1} considers the case that the driving L\'evy process $Y$ is spectrally positive, and Section~\ref{sec3.2} its spectrally negative counterpart.

\subsection{The spectrally positive case}
\label{sec3.1}
Suppose that $Y$ is spectrally positive, i.e., it has no downward jumps.
Define its Laplace exponent by $\varphi(\alpha):= \log {\mathbb E}\exp(-\alpha Y(1))$, and $\psi(\beta)$ its right-inverse; cf.\ \cite[Section 3.3]{Kyprianou}.
It is well known \cite[Section 6.5.2]{Kyprianou} that, for $\zeta\geqslant 0$,
\begin{equation}
{\mathbb E} \,{\rm e}^{-\alpha \overline{Y}(T_\zeta)} = \frac{\zeta}{\zeta - \varphi(\alpha)} \frac{\psi(\zeta) - \alpha}{\psi(\zeta)} .
\label{3.1}
\end{equation}
By substituting first $\zeta = \beta$ and then $\zeta = \beta + \omega$, we obtain from (\ref{3.1}) and Theorem~\ref{the2} the following expression for the Laplace-Stieltjes transform (LST) of $Y_{\beta,\omega}$.
\begin{proposition} \label{P1} If $Y$ is spectrally positive, then, for $\alpha\geqslant 0$,
\begin{equation}
{\mathbb E} \,{\rm e}^{-\alpha Y_{\beta,\omega}} = \frac{\alpha - \psi(\beta)}{\beta - \varphi(\alpha)} \frac{\beta}{\psi(\beta)} \frac{\beta + \omega - \varphi(\alpha)}{\alpha - \psi(\beta+\omega)} \frac{\psi(\beta+\omega)}{\beta+\omega}.
\label{3.2}
\end{equation}
\end{proposition}

\begin{remark}
\label{remark1}{\em 
We  also provide an alternative derivation of Proposition \ref{P1}, using a relation with an associated queueing model. First observe that $Z^-$ is minus the running minimum of $Y$ over an interval between two successive Poisson($\omega$) inspection epochs,
given that the latter epoch occurs {\em before} the killing epoch $T_\beta$; such an inspection interval is exp($\beta+\omega$) distributed. Hence $Z^-$ is exp($\psi(\beta+\omega)$) distributed, just like $-\underline{Y}(T_{\beta+\omega})$; see e.g.\ \cite[Section 6.5.2]{Kyprianou}.
Furthermore,
\begin{equation}{\mathbb E}\,e^{-\alpha Z} = {\mathbb E}\,e^{-\alpha Y(T_{\beta+\omega})}=\frac{\beta+\omega}{\beta+\omega - \varphi(\alpha)}.\end{equation}
Hence, because $Z = Z^+ - Z^-$, with $Z^+$ and $Z^-$ being independent:
\begin{equation}
{\mathbb E} \,{\rm e}^{-\alpha Z^+} = \frac{\beta+\omega}{\beta+\omega - \varphi(\alpha)} \frac{\psi(\beta+\omega) - \alpha}{\psi(\beta+\omega)} ;
\label{3.3}
\end{equation}
indeed, cf.\ (\ref{3.1}), $Z^+ \stackrel{{\rm d}}{=} \overline{Y}(T_{\beta+\omega})$, as we already noticed in (\ref{d6a}).
Now observe, cf.\ (\ref{2.1}), that  $Y_{\beta,\omega}$ can be interpreted as the waiting time of the $N_{\beta,\omega}$-th customer of an M/G/1 queue with generic interarrival time $Z^-$
and generic service time $Z^+$, with the first customer arriving in an empty system.
Its LST is given by
\begin{equation}\sum_{n=1}^{\infty} \left(\frac{\omega}{\beta+\omega}\right)^n \frac{\beta}{\beta+\omega} \,{\mathbb E} \,{\rm e}^{-\alpha W_n},\end{equation}
with $W_n$ denoting the waiting time of the $n$-th such customer. The next step is to use expression \cite[Eqn.\ (II.4.77)]{Cohen82} for the generating function of ${\mathbb E} \,{\rm e}^{-\alpha W_n}$.
After some elementary calculations, (\ref{3.2}) is recovered.} \hfill$\Diamond$
\end{remark}

\begin{remark}{\em 
Considering the special case that the spectrally positive L\'evy process is a compound Poisson process,
our model corresponds to the Cram\'er-Lundberg insurance risk model. 
In particular, taking $\varphi'(0)>0$ so that eventual  ruin is not certain, the LST of $Y_{\beta,\omega}$ with $\beta = 0$ immediately yields the LST of the (all-time) bankruptcy probability in the Cram\'er-Lundberg model with initial capital $u$.
That quantity was studied in \cite{AL} for exponentially distributed claim sizes, and in \cite{BEJ} for generally distributed claim sizes.

Through the duality relation between the Cram\'er-Lundberg model and its queueing counterpart, 
our results also provide insight into the M/G/1 queue. In particular, Theorem~\ref{the2} entails for the special case $\beta=0$ and $\varphi'(0)>0$ that 
the steady-state workload is distributed as the sum of two independent quantities, viz.\ (i) $\overline{Y}(T_\omega)$, the supremum of the workload  until the first Poisson($\omega$) inspection epoch, and (ii) 
$\overline{Y}_{0,\omega}$, the steady-state workload (or waiting time) in an M/G/1 queue with ${\rm exp}(\psi(\omega))$ distributed interarrival times and service times distributed as the $Z^+$ defined above (with $\beta=0$).} \hfill$\Diamond$
\end{remark}

\begin{remark}{\em 
In many real-life applications, it may be more natural to have inspection intervals that are constant, or at least have a small coefficient of variation, instead of being exponentially distributed.
In this remark we outline how one can use the alternative derivation of Proposition~\ref{P1}, as described in Remark~\ref{remark1}, to determine the LST of the running maximum at Erlang($k,k \omega$) distributed inspection moments, with $k\in{\mathbb N}$.
Note that the mean inspection interval still equals $1/\omega$, and that its squared coefficient of variation equals $1/k$, so that the case of a large $k$ emulates constant inspection intervals.
We shall again denote by $Y_{\beta,\omega}$ the running maximum at inspection moments until killing epoch $T_\beta \sim {\rm exp}(\beta)$,  and by $N_{\beta,\omega}$ the number of inspections before killing.
We now have
\begin{equation}
\label{3.3a}
{\mathbb P}(N_{\beta,\omega} = n) = \left(\frac{k \omega}{k \omega + \beta}\right)^{kn} \left(1 - \left(\frac{k \omega}{k \omega + \beta}\right)^k \right), \,\,\, n=0,1,\dots ~.
\end{equation}
Indeed, the first factor in the righthand side denotes the probability that at least $kn$ intervals $ \sim {\rm exp}(k \omega)$ occur before $T_\beta$, while the second factor denotes the probability that, subsequently, less than $k$ such exp$(k \omega)$ intervals occur, i.e., the ($n+1$)-th inspection does not occur before $T_\beta$.

Denote by $Z_1^-,\dots,Z_k^-$ the negatives of the running minima, and by $Z_1^+,\dots,Z_k^+$  the running maxima, over $k$ consecutive exp($k \omega$) intervals which together compose
one inspection interval.
Since $Y$ is a spectrally positive L\'evy process, it follows that $Z_1^-,\dots,Z_k^-$ are i.i.d.\ exp$(\psi(\beta + k \omega))$ distributed. We furthermore note that
$Z_1^+,\dots,Z_k^+$ are i.i.d., and the reasoning leading to (\ref{3.3}) shows that their LST is given by
\begin{equation}
{\mathbb E}\, {\rm e}^{-\alpha Z^+} = \frac{\beta+k  \omega}{\beta+ k \omega - \varphi(\alpha)} \frac{\psi(\beta+ k \omega) - \alpha}{\psi(\beta+ k \omega)} .
\label{3.3b}
\end{equation}
Moreover, all $Z_i^-$ and $Z_j^+$ are independent.
Now $Y_{\beta,\omega}$, as a supremum of partial sums $S_n = \sum_{m=1}^n Z_m$ with $Z_m = \sum_{i=1}^k Z_{m,i}^+ - \sum_{i=1}^k Z_{m,i}^-$, can be interpreted as the waiting time of the $N_{\beta,\omega}$-th customer of an E$_k$/G/1 queue with generic interarrival time $\sum_{i=1}^k Z_i^-$ and generic service time $\sum_{i=1}^k Z_i^+$, with 
the first customer arriving in an empty system. Its LST, and hence the LST of $Y_{\beta,\omega}$, is given by 
\begin{equation}
\sum_{n=1}^{\infty} {\mathbb P}(N_{\beta,\omega}=n) \,{\mathbb E} \,{\rm e}^{-\alpha W_n} \, = \,
\left(1 - \left(\frac{k \omega}{\beta+k \omega}\right)^k \right)
\sum_{n=1}^{\infty} \left(\frac{k \omega}{\beta+k \omega}\right)^{kn}
\,{\mathbb E} \,{\rm e}^{-\alpha W_n}.
\end{equation}
Finally, we can use  \cite[Thm.\ 25, p.\ 44]{Prabhu} for the generating function of ${\mathbb E} \,{\rm e}^{-\alpha W_n}$ in the E$_k$/G/1 queue.

We close this remark by once more focusing on the bankruptcy probability in insurance risk. In Subsection~\ref{sec2.2} we saw that the probability of bankruptcy before $T_\beta$ is given by
$\tilde{p}(u,T_\beta) = {\mathbb P}(Y_{\beta, \omega} > u)$. Hence  the LST of  the bankruptcy probability, in the case of a spectrally positive L\'evy process and Erlang($k,k \omega$) inspection intervals,
immediately follows from the LST of $Y_{\beta,\omega}$.
}
\hfill$\Diamond$
\end{remark}

\begin{remark}{\em 
The decomposition can be used to determine all moments of $Y_{\beta,\omega}$ from the (known) corresponding moments of $\overline Y(T_\beta)$ and $\overline Y(T_{\beta+\omega})$. In this remark we demonstrate this by providing such a computation for the mean and variance of $Y_{\beta,\omega}$.
Clearly, it suffices to be able to determine mean and variance of $\overline{Y}(T_{\zeta})$ for some $\zeta>0$:
\begin{equation}
{\mathbb E} \,\overline{Y}(T_\zeta) = \frac{1}{\psi(\zeta)} - \frac{\varphi'(0)}{\zeta} ,
\label{mean}
\end{equation}
\begin{equation}
{\mathbb V}{\rm ar}\, \overline{Y}(T_\zeta) = \frac{\varphi''(0)}{\zeta} + \left(\frac{\varphi'(0)}{\zeta}\right)^2 - \left(\frac{1}{\psi(\zeta)}\right)^2.
\label{vari}
\end{equation}
Due to the independence of the terms in the righthand side of the decomposition of Theorem~\ref{the2},
the mean and variance of $Y_{\beta,\omega}$ immediately follow by successively plugging in
$\zeta = \beta$ and $\zeta = \beta + \omega$ and subtracting the resulting expressions.} \hfill$\Diamond$
\end{remark}

\subsection{The spectrally negative case}
\label{sec3.2}
Suppose that $Y$ is spectrally negative, i.e., it has no upward jumps.
Consider an exp($\zeta$) distributed interval. Define the cumulant
$\Phi(\alpha) := \log {\mathbb E} \exp({\alpha Y(1)})$ and its right-inverse $\Psi(\beta).$
As follows directly from e.g.\ \cite[Section 6.5.2]{Kyprianou}, the running maximum $\overline Y (T_\zeta)$ is exponentially distributed 
with rate $\Psi(\zeta)$. 
Using this result with $\zeta = \beta$ and $\zeta = \beta + \omega$, and applying Theorem~\ref{the2}, we obtain an expression for the LST of $Y_{\beta,\omega}$ in the spectrally negative case.
\begin{proposition} \label{P2} If $Y$ is spectrally negative, then, for $\alpha\geqslant 0$,
\begin{equation}
{\mathbb E} \,{\rm e}^{-\alpha Y_{\beta,\omega}} = \frac{\Psi(\beta)}{\Psi(\beta) + \alpha} \frac{\Psi(\beta + \omega) + \alpha}{\Psi(\beta + \omega)} .
\label{3.4}
\end{equation}
\end{proposition}
Using Proposition \ref{P2}, an elementary computation reveals that $Y_{\beta,\omega}$ has an atom at zero, i.e., 
\begin{equation}{\mathbb P}(Y_{\beta,\omega} = 0) = \frac{\Psi(\beta)}{\Psi(\beta + \omega)},
\end{equation}
and is 
exp($\Psi(\beta)$) distributed with the complementary probability $1-{\mathbb P}(Y_{\beta,\omega} = 0)$ .

\begin{remark}{\em 
\label{remark2}
Just like in Remark~\ref{remark1}, we could also have obtained the LST of $Y_{\beta,\omega}$ by observing that $Y_{\beta,\omega}$ can be interpreted as the waiting time
of the $N_{\beta,\omega}$-th customer in a single-server queue with generic interarrival time $Z^-$ and generic service time $Z^+$, with the first customer arriving in an empty system.
In this case $Z^+$ is exp($\Psi(\beta + \omega)$) distributed. Hence we can now use \cite[Eqn.\ (II.3.100)]{Cohen82}
for the generating function of  ${\mathbb P}(W_n < s)$,
with $W_n$ denoting the waiting time of the $n$-th customer in the G/M/1 queue.
We close this remark by observing
that $Y_{0,\omega}$ is distributed as the steady-state waiting time  (if it exists) in the  above-described G/M/1 queue;
that waiting time also is exponentially distributed with an atom at zero.
} \hfill$\Diamond$
\end{remark}

\section{Asymptotics in the compound Poisson setting}
\label{sec4}
In this section we demonstrate the potential of our results by using them to establish the asymptotics of the bankruptcy probability for large initial capital $u$. The driving L\'evy process is a compound Poisson process with drift, characterized through its Laplace exponent
\begin{equation}
\label{4.1}
\varphi(\alpha) \, = \, r \alpha - \lambda(1-b(\alpha)) ;
\end{equation}
 here $\lambda$ can be viewed as the claim arrival rate, $r$ as the premium rate and $b(\cdot)$ as the LST of a generic claim size $B$.
We are interested in the behavior of $\tilde p(u)={\mathbb P}(Y_{0,\omega} > u)$ for large $u$, which can be interpreted as the bankruptcy probability  in the Cram\'er-Lundberg insurance risk model for
large initial capital $u$.
We assume that $\lambda \,{\mathbb E} B < r$, as otherwise ruin and bankruptcy are certain.
Section~\ref{sec4.1} treats the case of a light-tailed jump-size distribution, and Section~\ref{sec4.2} that of a heavy-tailed jump-size distribution.
Here light-tailed means that $b(\alpha)$ is finite for some $\alpha <0$, or equivalently, that ${\mathbb P}(B>x) = O({\rm e}^{-ax})$ for some $a>0$.
Heavy-tailed means that $b(\alpha)$ is infinite for all $\alpha <0$; we shall restrict ourselves to the well-known subclass ${\mathscr S}^\star$.

\subsection{The bankruptcy probability in the light-tailed case}
\label{sec4.1}
We assume in this subsection that $B$ is light-tailed in the sense that there is a unique strictly positive solution $\theta^\star$ of the equation $\varphi(-\theta^\star)=0$.
Our aim is to identify the asymptotic behavior of $\tilde p(u)$, and to compare it to the classic result for the asymptotic ruin probability $p(u)$ in the same  model (which coincides with the asymptotic behavior
of the waiting-time tail in the dual M/G/1 queue).
The well-known Cram\'er-Lundberg approximation  \cite[Theorem 5.3]{AA} states that
\begin{equation}
p(u) \, \sim \, \gamma {\rm e}^{-\theta^\star u} , 
\label{4.2}
\end{equation}
with $f(u)\sim g(u)$ denoting $f(u)/g(u)\to 1$ as $u\to\infty$,
and
\begin{equation}
\gamma \, = \, - \frac{\varphi'(0)}{\varphi'(-\theta^\star)} .
\label{4.3}
\end{equation}
To determine the asymptotics of $\tilde{p}(u)$ we return to an observation made in Remark~\ref{remark1}: $Y_{\beta,\omega}$ can be interpreted as the waiting time of the $N_{\beta,\omega}$-th customer of an M/G/1 queue with generic (exponentially distributed) interarrival time $Z^-$ and generic service time $Z^+$. For $\beta=0$, $Y_{0,\omega}$ becomes the steady-state waiting time for that queue.
Notice that ${\mathbb E}Z^- = 1/\psi(\omega)$ while, as $Z^+$ is distributed as $\overline{Y}(T_{\omega})$, it follows from (\ref{mean}) that ${\mathbb E}Z^+ = 1/\psi(\omega) - \varphi'(0)/\omega$.
As $\varphi'(0) = r-\lambda {\mathbb E} B > 0$, we have ${\mathbb E} Z^- > {\mathbb E} Z^+$, so the steady-state waiting-time distribution indeed exists.
As a consequence, 
we can  (again) rely on the Cram\'er-Lundberg approximation, or (equivalently) the tail asymptotics of the M/G/1 queue with generic interarrival time $Z^-$ and generic service time $Z^+$, cf.\ (\ref{4.2}):
\begin{equation}
\tilde{p}(u) \, \sim \, \tilde{\gamma} {\rm e}^{-\tilde{\theta}^\star u} , ~~~ u \rightarrow \infty .
\label{4.4}
\end{equation}
Our next task is to determine to determine the constants $\tilde{\theta}^\star$ and $\tilde{\gamma}$.
The customer arrival rate in the auxiliary M/G/1 queue is  $1/{\mathbb E}Z^- = \psi(\omega$),
the generic  service time $Z^+$ has  LST  (cf.\ (\ref{3.3}))
\[{\mathbb E} \,{\rm e}^{-\alpha Z^+} = \frac{\omega}{\omega - \varphi(\alpha)} \frac{\psi(\omega) - \alpha}{\psi(\omega)},\]
and the service speed or premium rate equals one.
Hence $\tilde{\theta}^\star$ is the unique positive solution of 
$\tilde{\varphi}(-\tilde{\theta}^\star)=0$, with 
$\tilde{\varphi}(\alpha) := \alpha - \psi(\omega)(1 - {\mathbb E}{\rm e}^{-\alpha Z^+})$:
\begin{equation} 
\-\tilde{\theta}^\star \, - \, \psi(\omega) \left(1 - \frac{\omega}{\psi(\omega)} \frac{\psi(\omega) + \tilde{\theta}^\star }{\varphi(-\tilde{\theta}^\star) - \omega}\right) \, = \, 0.
\label{4.5}
\end{equation}
It readily follows that $\tilde{\theta}^\star = \theta^\star$ satisfies this equation (as $\varphi(-\theta^\star)=0$), while otherwise there is only the negative solution $\tilde{\theta}^\star = - \psi(\omega)$.
The implication is that the bankruptcy probability $\tilde{p}(u)$ has the same decay rate $\theta^\star$ as the ruin probability $p(u)$.

Let us now determine the prefactor $\tilde{\gamma}$. Using (\ref{4.3}) with $\varphi(\alpha)$ replaced by
$\tilde{\varphi}(\alpha)$,
we find:
\begin{equation}
\tilde{\gamma} \, = \, - \frac{\tilde{\varphi}'(0)}{\tilde{\varphi}'(-\tilde{\theta}^\star)} \, = \, - \frac{1 - \psi(\omega) {\mathbb E}Z^+}{1 +\psi(\omega) \left.\frac{{\rm d}}{{\rm d} \alpha} {\mathbb E} \,{\rm e}^{-\alpha Z^+}\right|_{\alpha = -\tilde{\theta}^\star}} .
\label{4.6}
\end{equation}
A brief calculation, using (\ref{4.3}), results in
\begin{equation}
\tilde{\gamma} = \gamma \frac{\psi(\omega)}{\psi(\omega) + \theta^\star} .
\label{4.7}
\end{equation}
We have established the main result of this subsection.

\begin{proposition} Assume {$B$ is light-tailed.}
As $u\to\infty$,
\begin{equation}
\frac{\tilde{p}(u)}{p(u)}\to \gamma^\star_\omega:= \frac{\psi(\omega)}{\psi(\omega)+\theta^\star}.
\label{4.8}
\end{equation}
\end{proposition}

The above result shows that $\gamma^\star_\omega\uparrow 1$ as the inspection rate $\omega$ grows large. In addition, using that $\varphi(\alpha)-r\alpha +\lambda\to 0$ as $\alpha\to\infty$ implies that $\psi(\theta)- (\theta+\lambda)/r\to 0$ as $\theta\to\infty$, it can be verified that
\begin{equation}\lim_{\omega\to\infty}\omega(1-\gamma^\star_\omega) = r\theta^\star,\end{equation}
which means that, for $\omega$ large, $\gamma^\star_\omega$ behaves as $1- r\theta^\star/\omega$. This relation can be used to determine a `rule of thumb' by which one can determine the minimally required inspection rate $\omega$ such that the information loss due to Poisson inspection is below a given threshold. 

\subsection{The bankruptcy probability in the heavy-tailed case}
\label{sec4.2}
In this subsection we study the asymptotic behavior of the bankruptcy probability when $B$ is heavy-tailed. More specifically, we assume that $B  \in {\mathscr S}^\star$, a class
introduced by Kl\"uppelberg \cite{Kluppelberg}. A random variable $U$ on ${\mathbb R}$ belongs to ${\mathscr S}^\star$ iff its complementary distribution function $\overline F_U(x):={\mathbb P}(U>x)$ is positive for all $x$ and
\begin{equation}
\int_0^x \overline{F}_U(x-y) \overline{F}_U(y) \,{\rm d}y \, \sim \, 2 m_U\, \overline{F}_U(x)  ~~~ {\rm as} ~ x \rightarrow \infty ;
\label{4.9}
\end{equation}
here $m_U$ denotes the mean of $U$, restricted to the positive half-axis.
${\mathscr S}^\star$ is a class that is contained in, but is also very close to, the well-known class ${\mathscr S}$ of subexponential distributions.
The class ${\mathscr S}^\star$ has the convenient property that if $U \in {\mathscr S}^\star$ then both $U$ and $U^{{\rm res}}$, the latter random variable being characterized through
\begin{equation}
    {\mathbb P}(
U^{{\rm res}}\leqslant x) = \int_0^x \frac{{\mathbb P}(U>y)}{m_U} {\rm d}y,
\end{equation}
are subexponential.

In our analysis of the asymptotics of $\tilde p(u)$, we shall use a well-known result from \cite{EKM} concerning the supremum $M_{\sigma}$ of a random walk $\{S_n\}_{n\in{\mathbb N}}$ over an interval  $[0,\sigma]$, with $\sigma$
some random variable (see also \cite[Theorem 1]{FZ} for the more general case of $\sigma$ being a stopping time, and   \cite[p.\ 309]{AA} for the special case of a constant $\sigma$).
If the increments of the random walk are in ${\mathscr S}^\star$, with distribution function $F(\cdot)$ and complementary distribution function $\overline F(\cdot)$, then 
\begin{equation}
\lim_{u \rightarrow \infty} \frac{{\mathbb P}(M_{\sigma} >u)}{\overline{F}(u)} = {\mathbb E}\, \sigma .
\label{4.10}
\end{equation}

Notice that the increments attain values in ${\mathbb R}$, i.e., not necessarily in ${\mathbb R}^+$. If we take $F(\cdot)$ to be the distribution function of a claim size minus a claim interarrival time,
then  we can apply (\ref{4.10}) by considering the random number of claim arrivals in an exp($\omega$) interval.
This yields:
\begin{equation}
{\mathbb P}
\big(\overline{Y}(T_\omega) >u\big) \sim \frac{\lambda}{\omega} \bar{F}(u) \sim \frac{\lambda}{\omega}\, {\mathbb P}(B>u) .
\label{4.11}
\end{equation}
In order to use this result for determining the asymptotic behavior of the bankruptcy probability, we observe the following.
\begin{itemize}
    \item[(i)] $\overline{Y}(T_\omega)$ has the same distribution as $Z^+$, cf.\ (\ref{d6a}) with $\beta = 0$.
Hence $Z^+ \in {\mathscr S}^\star$, and $Z^{+,{\rm res}}$ is subexponential.
\item[(ii)] $Y_{0,\omega}$ can be viewed as the steady-state waiting time in an M/G/1 queue with generic service time $Z^+$ and arrival rate $\psi(\omega)$, 
so with load $\rho := {\mathbb E}Z^+/{\mathbb E}Z^-$.
Hence we can use a standard result for the waiting-time tail in the M/G/1 queue,
in which the residual service time 
$Z^{+,{\rm res}}$
is subexponential: \begin{equation}
  {\mathbb P}(Y_{0,\omega} >u) \sim \frac{\rho}{1-\rho} {\mathbb P}(Z^{+,{\rm res}}>u).  
\end{equation}
This result holds equivalently for the ruin probability in the dual Cram\'er-Lundberg model,
see e.g.\ \cite[Theorem X.2.1]{AA},
and hence we conclude:
\begin{align}
\tilde{p}(u) \, &= \,  {\mathbb P}(Y_{0,\omega} >u) \, \sim \, \frac{{\mathbb E} Z^+}{{\mathbb E}Z^- - {\mathbb E} Z^+} {\mathbb P}(Z^{+,{\rm res}} >u)
\nonumber
\\
&= \,
\frac{1}{{\mathbb E}Z^- - {\mathbb E}Z^+} \int_u^{\infty} {\mathbb P}(Z^+ >y)\, {\rm d} y
\nonumber
\\
& = \, \frac{\omega}{r - \lambda {\mathbb E}B} \int_u^{\infty} \frac{\lambda}{\omega}\, {\mathbb P}(B>y)\, {\rm d} y .
\label{4.12}
\end{align}
Here the last equality follows using ${\mathbb E}Z^- = 1/\psi(\omega)$ and (cf.\ (\ref{mean}))
${\mathbb E}Z^+ = 1/\psi(\omega) - \varphi'(0)/\omega$ with $\varphi'(0) = r-\lambda {\mathbb E}B >0$, and applying (\ref{4.11}).
\end{itemize}
Combining (i) and (ii) yields the main result of this subsection.

\begin{proposition} Assume $B\in {\mathscr S}^\star$.
As $u\to\infty$,
\begin{equation}
\tilde{p}(u) \, \sim \, \frac{\lambda {\mathbb E}B}{r - \lambda {\mathbb E}B} {\mathbb P}(B^{{\rm res}} >u) .
\label{4.13}
\end{equation}
\end{proposition}

Notice, again looking at  the statement in (ii) above regarding the tail behavior in an M/G/1 queue, that $\tilde p(u)$ has the exact same tail asymptotics as $p(u)$. In particular, the asymptotics of $\tilde p(u)$ do not depend on the inspection rate $\omega.$
This may look surprising at first sight, but realize that for $B\in {\mathscr S}^\star$ there is the intuition that  `a single big jump' is responsible (with overwhelming probability as $u\to\infty$) for exceeding a high level. This 
suggests that when ruin occurs, it is highly likely that the capital still is below zero at the next inspection moment.

\begin{remark}{\em 
In the special heavy-tailed case of regularly varying claim sizes we can identify the tail asymptotics of $\tilde{p}(u)$ in a more straightforward way by applying 
\cite[Theorem 8.1.6]{BGT} to the LST of $Y_{0,\omega}$ (as given by (\ref{3.4}) with $\beta=0$).}\hfill$\Diamond$
\end{remark}

\begin{remark}{\em 
As suggested in \cite{HJA},  
one could alternatively obtain the asymptotics of this section via  \cite[Eqn.\ (2)]{AI}, which in our notation reads as follows (see also (\ref{d1})):
\begin{equation}
\tilde{p}(u) = {\mathbb E} \, p(u+Z^-) .
\label{d1a}
\end{equation}
Use dominated convergence to take the limit inside the expectation of (\ref{d1a}) and
observe that $Z^-$ is exp$(\psi(\omega))$ distributed in this spectrally positive case.
In the light-tailed case one can then see that $\tilde{\theta}^\star = \theta^\star$ and that $\tilde{\gamma} = \gamma {\mathbb E} \,{ e}^{-\theta^\star Z^-}$, thus
yielding (\ref{4.7}).
In the heavy-tailed case the following holds for a subexponential random variable $X$: ${\mathbb P}(X>u) \sim {\mathbb P}(X>u+Z^-)$, and hence the subexponential asymptotics 
for the continuously observed and discretely observed processes are the same.
}\hfill$\Diamond$
\end{remark}

\section{Discussion and concluding remarks}
\label{sec5}
The main finding presented in this work is a decomposition result involving the running maximum value of a L\'evy process at Poisson inspection times as well as the running maximum of the permanently observed process. It allows the translation of known results for the permanently observed process in terms of the Poisson-observed process: such a translation procedure has been followed to find an explicit characterization of the running maximum value of spectrally one-sided L\'evy processes at Poisson inspection times, and to find tail asymptotics of the bankruptcy probability in the celebrated Cram\'er-Lundberg model.

Many extensions can be thought of. In the first place, one could think of the inter-inspection times {being phase-type}; to this end, potentially ideas from \cite[Section 5.1]{MOR} and \cite{AI2} can be used.
{It would also be interesting to see whether one can have similar decompositions as Theorem~\ref{the2} for overshoots, cf.\ \cite{AI}.
From a practical perspective, it would further be relevant to allow the inspection rate to depend on the process level.
Finally}, the connection with Poisson-observed L\'evy-driven storage systems can be explored, e.g.\ to study whether the results in the present paper can be used for hypothesis testing purposes \cite{MR}. 

\end{document}